\documentclass[12pt,thmsa]{article}%
\usepackage{amsmath}
\usepackage{amssymb}
\usepackage{sw20lart}
\usepackage{amsfonts}
\usepackage{graphicx}%
\begin{document}

\title{Geometric Algebras }
\author{{\footnotesize A. M. Moya}$^{2}$, {\footnotesize V. V. Fern\'{a}ndez}$^{1}%
${\footnotesize and W. A. Rodrigues Jr}.$^{1}${\footnotesize \ }\\$^{1}\hspace{-0.1cm}${\footnotesize Institute of Mathematics, Statistics and
Scientific Computation}\\{\footnotesize \ IMECC-UNICAMP CP 6065}\\{\footnotesize \ 13083-859 Campinas, SP, Brazil }\\{\footnotesize e-mail: walrod@ime.unicamp.br }\\{\footnotesize \ }$^{2}${\footnotesize Department of Mathematics, University
of Antofagasta, Antofagasta, Chile} \\{\footnotesize e-mail: mmoya@uantof.cl}}
\maketitle

\begin{abstract}
This is the first paper in a series of eight where in the first three we
develop a systematic approach to the geometric algebras of multivectors and
extensors, followed by five papers where those algebraic concepts are used in
a novel presentation of several topics of the differential geometry of
(smooth) manifolds of arbitrary global topology. A key tool for the
development of our program is the mastering of the euclidean geometrical
algebra of multivectors that is detailed in the present paper.

\end{abstract}
\tableofcontents

\section{Introduction}

This is the first paper in a following series of eight which have been
designed in order to show how Clifford (geometric) algebra methods can be
conveniently used in the study of differential geometry and geometrical
theories of the gravitational field. It dispenses the use of fiber bundle
theory\footnote{A reader interested in Clifford and spin-Clifford bundles may
consult, e.g., \cite{rodoliv2006}.} and is indeed a very powerful and economic
tool for performing sophisticated calculations. The first three
papers\footnote{This one and \cite{1,2}.} deal with the algebraic aspects of
the theory, namely Clifford algebras and the theory of extensors. Our
presentation is self  contained and serves besides the purpose of fixing our
conventions also the one of introducing a series of \ \textquotedblleft tricks
of the trade\textquotedblright\ (not found easily elsewhere) necessary for
quickly and efficient computations. The other five papers \cite{3,4,5,6,7}
develop a systematic approach to a theory of multivector and extensor calculus
and their use in the differential geometry of manifolds of arbitrary topology.
There are many novelties in our presentation, in particular the way we
introduce the concept of deformation of geometric structures, which is
discussed in detail in \cite{5} and which permit us to relate \cite{6,7} some
distinct geometric structures on a given manifold. Moreover, the method permit
also to solve problems in one given geometry in terms of an eventually simple
one, and here is a place where our theory may be a useful one for the study of
geometrical theories of the gravitational field.

The main issues discussed in the present paper are the constructions of the
euclidean and metric \textit{geometric} (or \textit{Clifford}) algebras of
\emph{multivectors} which can be associated to a real vector space $V$ of
dimension $n$ once we equip $V$ with an euclidean and an arbitrary non
degenerated metric of signature $(p,q)$ with $p+q=n$. The euclidean
geometrical algebra is the \textit{key} tool for performing almost all
calculations of the following papers. Metric geometric algebras are introduced
as \emph{deformations} of the euclidean geometric algebra. This is more
explored in \cite{2} where we introduce the concept of a
deformation\emph{\ extensor} associated to a given metric extensor and and by
proving the remarkable \emph{golden formula}. Extensors are a new kind of
geometrical objects which play a crucial role in the theory presented in this
series and the basics of their theory is described in \cite{2,3}. These
objects have been apparently introduced by Hestenes and Sobczyk in \cite{8}
and some applications of the concept appears in \cite{9}, but a rigorous
theory was lacking until recently\footnote{More details on the theory of
extensors may be found in \cite{10}.}. It is important to observe that in
\cite{8} the preliminaries of the geometric calculus have been applied to the
study of the differential geometry of \emph{vector} manifolds. However, as
admitted in \cite{11} there are some problems with that approach. In contrast,
our formulation applies to manifolds of a arbitrary topology and is free of
the problems that paved the construction in \cite{8}. There are many novelties
and surprises in what follows, e.g., the concept of \textit{deformed
geometries} relative to a given geometry, an \emph{intrinsic} Cartan calculus
and other topics that are ready to be used in geometrical theories of the
gravitational field, as the reader will convince himself consulting the other
papers in the series. 

As for the explicit contents of the present paper we introduce the concept of
\textit{multivectors} in Section 2 and their \textit{exterior algebra} in
section 3. In Section 4 we introduce the \textit{scalar product} of
multivectors and in section 5 the concepts of right and left
\textit{contractions} and \textit{interior} algebras. In Section 6 we give a
definition of a general \emph{real} Clifford (or geometrical) algebra of
multivectors and in Section 7 we study in details the relation between the
euclidean and pseudo-euclidean geometrical algebras. In section 8 we present
our conclusions. As additional references to several aspects of the theory of
the theory of Clifford algebras, that eventually may help the interested
reader, we quote\footnote{For a presentation of the theory of algebraic and
Dirac-Hestenes spinors and spinor fields, we quote \cite{rod041,moro}.
Applications of the theory may be found in \cite{rodoliv2006} and in
references quoted there.} \cite{10,12,13,14,15,rod041,rodoliv2006}.

\section{Multivectors}

Let $\bigwedge V$ denote the Cartesian product of all the $k$-vector spaces
$\bigwedge^{k}V,$ with $0\leq k\leq n,$ i.e., $\bigwedge V=\mathbb{R}\times
V\times\ldots\times\bigwedge^{n}V$. Here $V$ is $n$-dimensional vector space
over the real field $\mathbb{R}$.%
\begin{equation}
\bigwedge V=\{(X_{0},X_{1},\ldots,X_{n})\text{ }|\text{ }X_{k}\in\bigwedge
^{k}V,\text{ for each }k=0,1,\ldots,n\}. \label{M.0}%
\end{equation}

We introduce an equality among the $(n+1)$-uples of $\bigwedge V,$
\begin{equation}
(X_{0},X_{1},\ldots,X_{n})=(Y_{0},Y_{1},\ldots,Y_{n})\Leftrightarrow
X_{0}=Y_{0},\text{ }X_{1}=Y_{1},\ldots,\text{ }X_{n}=Y_{n}. \label{M.1}%
\end{equation}

$\bigwedge V$ has a structure of a real vector space, naturally induced by the
real vector space structure of all the $k$-vector spaces $\bigwedge^{k}V$. It
is realized by defining:

(\textbf{i}) The addition of $(X_{0},X_{1},\ldots,X_{n})\in\bigwedge V$ and
$(Y_{0},Y_{1},\ldots,Y_{n})\in\bigwedge V:$%
\begin{equation}
(X_{0},X_{1},\ldots,X_{n})+(Y_{0},Y_{1},\ldots,Y_{n})=(X_{0}+Y_{0},X_{1}%
+Y_{1}\ldots,X_{n}+Y_{n})\in\bigwedge V. \label{M.2a}%
\end{equation}

(\textbf{ii}) The scalar multiplication of $(X_{0},X_{1},\ldots,X_{n}%
)\in\bigwedge V$ by $\alpha\in\mathbb{R}:$%
\begin{equation}
\alpha(X_{0},X_{1},\ldots,X_{n})=(\alpha X_{0},\alpha X_{1},\ldots,\alpha
X_{n})\in\bigwedge V. \label{M.2b}%
\end{equation}

The vectors in $\bigwedge V$ will be called \emph{multivectors over} $V$.
Sometimes, they will be named as \emph{multivectors of} $\bigwedge V.$ The
real vector space $\bigwedge V$ is called the \emph{space of multivectors
over} $V$.

The zero multivector for $\bigwedge V$ is just $0=(0_{0},0_{1},\ldots,0_{n}),$
where $0_{k}$ denotes the zero $k$-vector of $\bigwedge^{k}V$.

For each $k=0,1,\ldots,n$ the linear mapping $\pi_{k}:\bigwedge V\rightarrow
\bigwedge^{k}V$ such that
\begin{equation}
\text{if }X=(X_{0},X_{1},\ldots,X_{n}),\text{ then }\pi_{k}(X)=X_{k}
\label{M.3}%
\end{equation}
is called the $k$\emph{-component projection operator} because $\pi_{k}(X)$ is
just the $k$-component of $X$.

For each $k=0,1,\ldots,n$ any multivector $X$ such that $\pi_{j}%
(X_{(k)})=0_{j},$ for $k\neq j,$ is said to be an \emph{homogeneous
multivector of degree }$k,$ or for short, a $k$\emph{-homogeneous
multivector}. That means that any component of $X$ which is not its
$k$-component is necessarily zero, but the $k$-component of $X$ may be zero or
not. Thus, any $k$-homogeneous multivector might be denoted by
\begin{equation}
X_{(k)}=(\ldots,0_{k-1},X_{k},0_{k+1},\ldots), \label{M.4}%
\end{equation}
where $X_{k}$ is some $k$-vector belonging to $\bigwedge^{k}V$.

It should be noticed that $0$ is an homogeneous multivector of any degree
$0,1,\ldots,n.$

The set of homogeneous multivectors of degree $k,$ i.e.,
\begin{equation}
\bigwedge\nolimits^{(k)}V=\{X_{(k)}=(\ldots,0_{k-1},X_{k},0_{k+1}%
,\ldots)\text{ }|\text{ }X_{k}\in\bigwedge\nolimits^{k}V\}\subset\bigwedge V
\label{M.5}%
\end{equation}
is a vector subspace of $\bigwedge V.$ Indeed, we have that $0\in
\bigwedge^{(k)}V$ and the addition of $k$-homogeneous multivectors and the
scalar multiplications of $k$-homogeneous multivectors by real numbers are
themselves $k$-homogeneous multivectors.

It is easy to see that $\bigwedge^{(k)}V$ is linearly isomorphic
$\bigwedge^{k}V.$ For each $k=0,1,\ldots,n$ there exists a linear isomorphism
between $\bigwedge^{k}V$ and $\bigwedge^{(k)}V$ which is realized by defining
$\tau_{k}:\bigwedge^{k}V\rightarrow\bigwedge^{(k)}V$ and $\tau_{k}%
^{-1}:\bigwedge^{(k)}V\rightarrow\bigwedge^{k}V$ such that
\begin{align}
\tau_{k}(X_{k})  &  =X_{(k)},\label{M.6a}\\
\tau_{k}^{-1}(X_{(k)})  &  =X_{k}. \label{M.6b}%
\end{align}
We have indeed that both of $\tau_{k}$ and $\tau_{k}^{-1}$ are linear
mappings, and they are inverses to each other, i.e., $\tau_{k}^{-1}\circ
\tau_{k}=i_{\bigwedge^{k}V}$ and $\tau_{k}\circ\tau_{k}^{-1}=i\bigwedge
^{(k)}V,$ where $i_{\bigwedge^{k}V}$ and $i_{\bigwedge^{(k)}V}$ are the
respective identity functions for $\bigwedge^{k}V$ and $\bigwedge^{(k)}V.$
Note that $\tau_{k}(\bigwedge^{k}V)=\bigwedge\limits^{(k)}V$ and
$\bigwedge^{k}V=\tau_{k}^{-1}(\bigwedge^{(k)}V).$

Such a linear isomorphism $\tau_{k},$ named as the $k$-isomorphism for short
in what follows, is a key piece in a our theory of the geometric algebra of multivectors.

Then,
\begin{equation}
\dim\bigwedge\nolimits^{(k)}V=\dim\bigwedge\nolimits^{k}V=\binom{n}{k}
\label{M.7}%
\end{equation}

Any $X=(X_{0},X_{1},\ldots,X_{n})\in\bigwedge V$ can be written as sum of all
their $k$-homogeneous multivectors $\tau_{k}(X_{k})=\tau_{k}\circ\pi_{k}%
(X)\in\bigwedge^{(k)}V,$ i.e.,
\begin{equation}
X=\overset{n}{\underset{k=0}{\sum}}\tau_{k}(X_{k})=\overset{n}{\underset
{k=0}{\sum}}\tau_{k}\circ\pi_{k}(X). \label{M.8}%
\end{equation}
To prove Eq.(\ref{M.8}) we should use Eqs.(\ref{M.2a}), (\ref{M.4}) and
(\ref{M.6a}), and Eq.(\ref{M.3}).

If there exists $X\in\bigwedge V$ such that $\pi_{k}(X)=0_{k}$ for $k\neq l,$
and $\pi_{l}(X)=X_{l},$ then
\begin{equation}
X=\tau_{l}(X_{l}). \label{M.8a}%
\end{equation}
It is an immediate consequence of \texttt{Eq.}(\ref{M.8}).

Eq.(\ref{M.8}) implies that the space $\bigwedge V$ can be written as
sum\footnote{Recall that if $S_{1}$ and $S_{2}$ are subspaces of any space
$W,$ then $S_{1}+S_{2}$ is just the subspace of $W$ defined by $S_{1}%
+S_{2}=\{v_{1}+v_{2}$ $|$ $v_{1}\in S_{1}$ and $v_{2}\in S_{2}\}.$} of all
their $k$-subspaces $\bigwedge^{(k)}V,$ i.e.,
\begin{equation}
\bigwedge V=\bigwedge\nolimits^{(0)}V+\bigwedge\nolimits^{(1)}V+\ldots
+\bigwedge\nolimits^{(n)}V, \label{M.9a}%
\end{equation}
or into a more suggestive form,
\begin{equation}
\bigwedge V=\overset{n}{\underset{k=0}{\sum}}\tau_{k}(\bigwedge\nolimits^{(k)}%
V)=\overset{n}{\underset{k=0}{\sum}}\tau_{k}\circ\pi_{k}(\bigwedge V).
\label{M.9aa}%
\end{equation}
But, since
\begin{equation}
\bigwedge\nolimits^{(j)}V\cap\bigwedge\nolimits^{(k)}V=\{0\},\text{ for }j\neq
k, \label{M.9b}%
\end{equation}
we see that Eq\texttt{.}(\ref{M.9a}) can still be written as
\begin{equation}
\bigwedge V=\bigwedge\nolimits^{(0)}V\oplus\bigwedge\nolimits^{(1)}%
V\oplus\ldots\oplus\bigwedge\nolimits^{(n)}V. \label{M.10}%
\end{equation}

As an immediate consequence\footnote{Recall that $\dim(S_{1}+S_{2})=\dim
S_{1}+\dim S_{2}-\dim(S_{1}\bigcap S_{2}).$} of Eqs.(\ref{M.7}), (\ref{M.9a})
and (\ref{M.9b}) we have that
\begin{equation}
\dim\bigwedge V=\binom{n}{0}+\binom{n}{1}+\ldots+\binom{n}{n}=2^{n}.
\label{M.11}%
\end{equation}

For each $k=0,1,\ldots,n$ the linear mapping $\left\langle \left.  {}\right.
\right\rangle _{k}:\bigwedge V\rightarrow\bigwedge V$ such that
\begin{equation}
\text{if }X=(X_{0},X_{1},\ldots,X_{n}),\text{ then }\left\langle
X\right\rangle _{k}=X_{(k)}, \label{M.12}%
\end{equation}
is called the $k$-\emph{part operator}. $\left\langle X\right\rangle _{k}$ is
read as the $k$-\emph{part of} $X.$

The $k$-component projection operator, the $k$-\emph{isomorphism} and the
$k$-part operator are involved in the following basic properties
\begin{equation}
\pi_{j}\circ\tau_{k}(X_{k})=\left\{
\begin{array}
[c]{cc}%
0_{j}, & \text{if }j\neq k\\
X_{k}, & \text{if }j=k
\end{array}
\right.  \label{M.15a}%
\end{equation}
and
\begin{equation}
\tau_{k}\circ\pi_{k}=\left\langle \left.  {}\right.  \right\rangle _{k},\text{
for each }k=0,1,\ldots,n. \label{M.15b}%
\end{equation}

We notice that \texttt{Eq.}(\ref{M.15a}) implies that
\begin{equation}
\pi_{k}\circ\tau_{k}=i_{\bigwedge^{k}V},\text{ for each }k=0,1,\ldots,n,
\label{M.16}%
\end{equation}
i.e., $\pi_{k}$ is the \emph{left inverse} of $\tau_{k}$.

In order to prove \texttt{Eq.}(\ref{M.15a}) we should use Eq.(\ref{M.6a}) and
take into account the definition of $k$-homogeneous multivectors, as set by
Eq.(\ref{M.4}).

Eq.(\ref{M.15b}) follows directly of using Eqs.(\ref{M.3}), (\ref{M.6a}) and
(\ref{M.12}).

Note that $X\in\bigwedge V$ is an homogeneous multivector of degree $k,$ i.e.,
$X\in\bigwedge^{(k)}V,$ if and only if
\begin{equation}
X=\left\langle X\right\rangle _{k}. \label{M.17}%
\end{equation}
This logical equivalence is an immediate consequence of Eq.(\ref{M.15b})
whenever Eq.(\ref{M.3}) and Eq.(\ref{M.6a}) are taken into account.

Any multivector can be written as sum of all their own $k$-parts, i.e.,
\begin{equation}
X=\overset{n}{\underset{k=0}{\sum}}\left\langle X\right\rangle _{k}.
\label{M.18}%
\end{equation}
This result follows directly from Eq.(\ref{M.8}) and Eq.(\ref{M.15b}).

\section{Exterior Product}

Let $TV=\sum\limits_{k=0}^{\infty}T^{k}V$ be the tensor algebra of $V$ and
$\mathcal{A}$ the so-called \emph{antisymmetrization operator}, i.e., a linear
mapping $\mathcal{A}:T^{k}V\rightarrow\bigwedge^{k}V$ such that

(\textbf{i}) for all $\alpha\in\mathbb{R}:$%
\begin{equation}
\mathcal{A}\alpha=\alpha, \label{EP.2a}%
\end{equation}

(\textbf{ii}) for all $v\in V:$%
\begin{equation}
\mathcal{A}v=v, \label{EP.2b}%
\end{equation}

(\textbf{ii}) for all $t\in T^{k}V,$ with $k\geq2:$%
\begin{equation}
\mathcal{A}t(\omega^{1},\ldots,\omega^{k})=\frac{1}{k!}\epsilon_{i_{1}\ldots
i_{k}}t(\omega^{i_{1}},\ldots,\omega^{i_{k}}), \label{EP.2c}%
\end{equation}
where $\epsilon_{i_{1}\ldots i_{k}}$ is the so-called \emph{permutation symbol
of order} $k,$
\begin{equation}
\epsilon_{i_{1}\ldots i_{k}}\equiv\epsilon^{i_{1}\ldots i_{k}}=\left\{
\begin{array}
[c]{cc}%
1, & \text{if }i_{1}\ldots i_{k}\text{ is even permutation of }1\ldots k\\
-1, & \text{if }i_{1}\ldots i_{k}\text{ is odd permutation of }1\ldots k\\
0, & \text{otherwise}%
\end{array}
\right.  . \label{EP.2d}%
\end{equation}
The exterior product of $X_{p}\in\bigwedge^{p}V$ and $Y_{q}\in\bigwedge^{q}V,$
namely $X_{p}\wedge Y_{q}\in\bigwedge^{p+q}V,$ is defined\footnote{Take
noticie that other definitions with other factors before the
antisymmentrization operator are possible.} by
\begin{equation}
X_{p}\wedge Y_{q}=\frac{(p+q)!}{p!q!}\mathcal{A}(X_{p}\otimes Y_{q}),
\label{EP.1}%
\end{equation}
where $X_{p}\otimes Y_{q}$ is the tensor product of $X_{p}$ and $Y_{q},$

Let $\{e_{j}\}$ be a basis of $V,$ and $\{\varepsilon^{j}\}$ be its dual basis
for $V^{*},$ i.e., $\varepsilon^{j}(e_{i})=\delta_{i}^{j}.$ Now, let us take
$t\in T^{k}V$ with $k\geq1.$ Such a \emph{contravariant }$k$-\emph{tensor} $t$
can be expanded onto the $k$-\emph{tensor basis} $\{e_{j_{1}}\otimes
\ldots\otimes e_{j_{k}}\}$ with $j_{1},\ldots,j_{k}=1,\ldots,n$ by the
well-known formula
\begin{equation}
t=t^{j_{1}\ldots j_{k}}e_{j_{1}}\otimes\ldots\otimes e_{j_{k}}, \label{EP.2e}%
\end{equation}
where $t^{j_{1}\ldots j_{k}}=t(\varepsilon^{j_{1}},\ldots,\varepsilon^{j_{k}%
})$ are the so-called $j_{1}\ldots j_{k}$-\emph{contravariant components} of
$t$ with respect to $\{e_{j_{1}}\otimes\ldots\otimes e_{j_{k}}\}.$

From Eq.(\ref{EP.2c}) it follows (non-trivially) a remarkable identity which
holds for the basis $1$-forms $\varepsilon^{1},\ldots,\varepsilon^{n}$
belonging to $\{\varepsilon^{j}\}.$ It is
\begin{equation}
\mathcal{A}t(\varepsilon^{j_{1}},\ldots,\varepsilon^{j_{k}})=\frac{1}%
{k!}\delta_{i_{1}\ldots i_{k}}^{j_{1}\ldots j_{k}}t(\varepsilon^{i_{1}}%
,\ldots,\varepsilon^{i_{k}}), \label{EP.2f}%
\end{equation}
where $\delta_{i_{1}\ldots i_{k}}^{j_{1}\ldots j_{k}}$ is the so-called
\emph{generalized Kronecker symbol of order} $k,$
\begin{equation}
\delta_{i_{1}\ldots i_{k}}^{j_{1}\ldots j_{k}}=\det\left[
\begin{array}
[c]{ccc}%
\delta_{i_{1}}^{j_{1}} & \ldots & \delta_{i_{1}}^{j_{k}}\\
\ldots & \ldots & \ldots\\
\delta_{i_{k}}^{j_{1}} & \ldots & \delta_{i_{k}}^{j_{k}}%
\end{array}
\right]  \text{ with }i_{1},\ldots,i_{k}\text{ and }j_{1},\ldots,j_{k}\text{
running from }1\text{ to }n. \label{EP.2g}%
\end{equation}

Let us take $X\in\bigwedge^{k}V$ with $k\geq2.$ By definition $X\in T^{k}V$
and is completely skew-symmetric, hence, it must be $X=\mathcal{A}X.$ Then, by
using Eq.(\ref{EP.2f}) we get a combinatorial identity which relates the
$i_{1}\ldots i_{k}$-components to the $j_{1}\ldots j_{k}$-components for $X.$
It is
\begin{equation}
X^{j_{1}\ldots j_{k}}=\frac{1}{k!}\delta_{i_{1}\ldots i_{k}}^{j_{1}\ldots
j_{k}}X^{i_{1}\ldots i_{k}}. \label{EP.2h}%
\end{equation}

From Eq.(\ref{EP.2c}) by using a well-known property of the antisymmetrization
operator, namely: $\mathcal{A}(\mathcal{A}t\otimes u)=\mathcal{A}%
(t\otimes\mathcal{A}u)=\mathcal{A}(t\otimes u),$ a noticeable formula for
expressing simple $k$-vectors in terms of the tensor products of $k$ vectors
can be easily deduced. It is
\begin{equation}
v_{1}\wedge\ldots\wedge v_{k}=\epsilon^{i_{1}\ldots i_{k}}v_{i_{1}}%
\otimes\ldots\otimes v_{i_{k}}. \label{EP.3}%
\end{equation}
If $\omega^{1},\ldots,\omega^{k}\in V^{\ast},$ then
\begin{equation}
v_{1}\wedge\ldots\wedge v_{k}(\omega^{1},\ldots,\omega^{k})=\epsilon
^{i_{1}\ldots i_{k}}\omega^{1}(v_{i_{1}})\ldots\omega^{k}(v_{i_{k}}).
\label{EP.4}%
\end{equation}

Eq.(\ref{EP.3}) implies (non-trivially) a remarkable identity which holds for
the basis vectors $e_{1},\ldots,e_{n}$ belonging to any basis $\{e_{j}\} $ of
$V.$ It is
\begin{equation}
e_{i_{1}}\wedge\ldots\wedge e_{i_{k}}=\delta_{i_{1}\ldots i_{k}}^{j_{1}\ldots
j_{k}}e_{j_{1}}\otimes\ldots\otimes e_{j_{k}}. \label{EP.5a}%
\end{equation}

Once again let us take $X\in\bigwedge^{k}V$ with $k\geq2.$ Since $X\in T^{k}V$
and is completely skew-symmetric, the use of Eq.(\ref{EP.2h}) and
Eq.(\ref{EP.5a}) in Eq.(\ref{EP.2e}) allows us to obtain the expansion
formula
\begin{equation}
X=\frac{1}{k!}X^{i_{1}\ldots i_{k}}e_{i_{1}}\wedge\ldots\wedge e_{i_{k}}.
\label{EP.5b}%
\end{equation}

We recall now the basic properties of the elementary exterior product of
$p$-vector and $q$-vector.

For any $V_{p},W_{p}\in\bigwedge^{p}V$ and $X_{q},Y_{q}\in\bigwedge^{q}V$
\begin{align}
(V_{p}+W_{p})\wedge X_{q}  &  =V_{p}\wedge X_{q}+W_{p}\wedge X_{q}%
,\label{EP.5c}\\
V_{p}\wedge(X_{q}+Y_{q})  &  =V_{p}\wedge X_{q}+V_{p}\wedge Y_{q}\text{
(distributive laws).} \label{EP.5d}%
\end{align}

For any $X_{p}\in\bigwedge^{p}V,$ $Y_{q}\in\bigwedge^{q}V$ and $Z_{r}%
\in\bigwedge^{r}V$
\begin{equation}
(X_{p}\wedge Y_{q})\wedge Z_{r}=X_{p}\wedge(Y_{q}\wedge Z_{r})\text{
(associative law).} \label{EP.5e}%
\end{equation}

For any $X_{p}\in\bigwedge^{p}V$ and $Y_{q}\in\bigwedge^{q}V$
\begin{equation}
X_{p}\wedge Y_{q}=(-1)^{pq}Y_{q}\wedge X_{p}. \label{EP.5f}%
\end{equation}

\subsection{Exterior Product of Multivectors}

The exterior product of $X,Y\in\bigwedge V,$ namely $X\wedge Y\in\bigwedge V,$
is defined by
\begin{equation}
X\wedge Y=\overset{n}{\underset{k=0}{\sum}}\overset{k}{\underset{j=0}{\sum}%
}\tau_{k}(\pi_{j}(X)\wedge\pi_{k-j}(Y)). \label{EP.6}%
\end{equation}
Note that on the right side there appears the exterior product of $j$-vectors
and $(k-j)$-vectors, as defined by Eq.(\ref{EP.1}), which means that
\begin{equation}
\pi_{k}(X\wedge Y)=\overset{k}{\underset{j=0}{\sum}}\pi_{j}(X)\wedge\pi
_{k-j}(Y),\text{ for each }k=0,1,\ldots,n \label{EP.6a}%
\end{equation}
and, if $X=(X_{0},\ldots,X_{k},\ldots,X_{n})$ and $Y=(Y_{0},\ldots
,Y_{k},\ldots,Y_{n}),$ then
\begin{equation}
X\wedge Y=(X_{0}Y_{0},\ldots,\overset{k}{\underset{j=0}{\sum}}X_{j}\wedge
Y_{k-j},\ldots,\overset{n}{\underset{j=0}{\sum}}X_{j}\wedge Y_{n-j}).
\label{EP.6b}%
\end{equation}

This exterior product is an internal law on $\bigwedge V.$ It is associative
and satisfies the usual distributive laws (on the left and on the right).

In order to prove that the exterior product of multivectors, as defined by
Eq.(\ref{EP.6}), satisfies the associative law we should use the summation
identity $\overset{k}{\underset{j=0}{%
{\displaystyle\sum}
}}\overset{j}{\underset{i=0}{%
{\displaystyle\sum}
}}X_{i}\wedge(Y_{j-i}\wedge Z_{k-j})=\overset{k}{\underset{j=0}{%
{\displaystyle\sum}
}}\overset{k-j}{\underset{i=0}{%
{\displaystyle\sum}
}}X_{j}\wedge(Y_{i}\wedge Z_{k-j-i})$ and Eq.(\ref{EP.5e}). Thus, by using
Eq.(\ref{EP.6a}) a straightforward calculation gives
\begin{align*}
\pi_{k}((X\wedge Y)\wedge Z)  &  =\overset{k}{\underset{j=0}{%
{\displaystyle\sum}
}}\pi_{j}(X\wedge Y)\wedge Z_{k-j}=\overset{k}{\underset{j=0}{%
{\displaystyle\sum}
}}(\overset{j}{\underset{i=0}{%
{\displaystyle\sum}
}}X_{i}\wedge Y_{j-i})\wedge Z_{k-j}\\
&  =\overset{k}{\underset{j=0}{%
{\displaystyle\sum}
}}\overset{j}{\underset{i=0}{%
{\displaystyle\sum}
}}(X_{i}\wedge Y_{j-i})\wedge Z_{k-j}=\overset{k}{\underset{j=0}{%
{\displaystyle\sum}
}}\overset{j}{\underset{i=0}{%
{\displaystyle\sum}
}}X_{i}\wedge(Y_{j-i}\wedge Z_{k-j})\\
&  =\overset{k}{\underset{j=0}{%
{\displaystyle\sum}
}}\overset{k-j}{\underset{i=0}{%
{\displaystyle\sum}
}}X_{j}\wedge(Y_{i}\wedge Z_{k-j-i})=\overset{k}{\underset{j=0}{%
{\displaystyle\sum}
}}X_{j}\wedge(\overset{k-j}{\underset{i=0}{%
{\displaystyle\sum}
}}Y_{i}\wedge Z_{k-j-i})\\
&  =\overset{k}{\underset{j=0}{%
{\displaystyle\sum}
}}X_{j}\wedge\pi_{k-j}(Y\wedge Z)=\pi_{k}(X\wedge(Y\wedge Z)).
\end{align*}

The distributive laws on the left and on the right are immediate consequences
of the Eq.(\ref{EP.5c}) and Eq.(\ref{EP.5d}), respectively.

The space of multivectors $\bigwedge V$ endowed with this exterior product
$\wedge$ is an associative algebra called the \emph{exterior algebra of
multivectors}.

Let us take $X_{p}\in\bigwedge^{p}V$ and $Y_{q}\in\bigwedge^{q}V$ with $0\leq
p+q\leq n.$ By using Eq.(\ref{EP.6a}) and Eq.(\ref{M.15a}), we have that
\begin{align*}
\pi_{k}(\tau_{p}(X_{p})\wedge\tau_{q}(Y_{q}))  &  =\overset{k}{\underset
{j=0}{\sum}}\pi_{j}\circ\tau_{p}(X_{p})\wedge\pi_{k-j}\circ\tau_{q}(Y_{q})\\
&  =\overset{k}{\underset{j=0}{\sum}}\left\{
\begin{array}
[c]{cc}%
0_{j}, & j\neq p\\
X_{p}, & j=p
\end{array}
\right.  \wedge\left\{
\begin{array}
[c]{cc}%
0_{k-j}, & k-j\neq q\\
Y_{q}, & k-j=q
\end{array}
\right.  .
\end{align*}

But, the sum on the right side of the equation above is $0_{k}=\overset
{k}{\underset{j=0}{\sum}}0_{j}\wedge0_{k-j},$ for each $k=0,1,\ldots,n$ unless
there exists some $k_{0}$ with $0\leq k_{0}\leq n$ such that $k_{0}-p=q.$
Whence, we see that $k_{0}=p+q$ is the unique number which can satisfy the
required conditions. Then, from the above equation it follows that
\begin{align*}
\pi_{k}(\tau_{p}(X_{p})\wedge\tau_{q}(Y_{q}))  &  =0_{k},\text{ for }k\neq
p+q\\
\pi_{p+q}(\tau_{p}(X_{p})\wedge\tau_{q}(Y_{q}))  &  =X_{p}\wedge Y_{q}.
\end{align*}
Hence, using Eq.(\ref{M.8a}) we finally get
\begin{equation}
\tau_{p}(X_{p})\wedge\tau_{q}(Y_{q})=\tau_{p+q}(X_{p}\wedge Y_{q}).
\label{EP.7}%
\end{equation}

\section{Metric Structure}

Let us equip $V$ with a \emph{metric tensor}, i.e., a symmetric and
non-degenerate covariant $2$-tensor over $V,$ $G:V\times V\rightarrow
\mathbb{R}$ such that
\begin{align}
G(v,w)  &  =G(w,v)\text{ for all }v,w\in V.\label{MS.1a}\\
\text{If }G(v,w)  &  =0\text{ for all }w\in V,\text{ then }v=0. \label{MS.1b}%
\end{align}

As usual we write
\begin{equation}
G(v,w)\equiv v\cdot w, \label{MS.1c}%
\end{equation}
and call $v\cdot w$ the \emph{scalar product of the vectors }$v,w\in V.$

The pair $(V,G)$ is called a \emph{metric structure} for $V.$ Sometimes, $V$
is said to be a \emph{scalar product vector space.}

Let $\{e_{k}\}$ be any basis of $V,$ and $\{\varepsilon^{k}\}$ be its dual
basis for $V^{*}.$ As we know, $\{\varepsilon^{k}\}$ is the unique basis of
$V^{*}$ which satisfies $\varepsilon^{k}(e_{j})=\delta_{j}^{k}.$

Let $G_{jk}=G(e_{j},e_{k}),$ since $G$ is non-degenerate, it follows that
$\det\left[  G_{jk}\right]  \neq0.$ Then, there exist the $jk$-entries for the
inverse matrix of $\left[  G_{jk}\right]  ,$ namely $G^{jk},$ i.e.,
$G^{ks}G_{sj}=G_{js}G^{sk}=\delta_{j}^{k}.$

We introduce the \emph{scalar product} \emph{of} $1$-\emph{forms}
$\omega,\sigma\in V^{\ast}$ by
\begin{equation}
\omega\cdot\sigma=G^{jk}\omega(e_{j})\sigma(e_{k}). \label{MS.2}%
\end{equation}
It should be noticed that the real number given by Eq.(\ref{MS.2}) does not
depend on the choice of $\{e_{k}\}$.

Now, we can define the so-called reciprocal bases of $\{e_{k}\}$ and
$\{\varepsilon^{k}\}.$ Associated to $\{e_{k}\}$ we introduce the well-defined
basis $\{e^{k}\}$ by
\begin{equation}
e^{k}=G^{ks}e_{s},\text{ for each }k=1,\ldots,n. \label{MS.3a}%
\end{equation}
Such $e^{1},\ldots,e^{n}\in V$ are the unique basis vectors for $V$ which
satisfy
\begin{equation}
e^{k}\cdot e_{j}=\delta_{j}^{k}. \label{MS.3b}%
\end{equation}
Associated to $\{\varepsilon^{k}\},$ we can also introduce a well-defined
basis $\{\varepsilon_{k}\}$ by
\begin{equation}
\varepsilon_{k}=G_{ks}\varepsilon^{s},\text{ for each }k=1,\ldots,n.
\label{MS.4a}%
\end{equation}
Such $\varepsilon_{1},\ldots,\varepsilon_{n}\in V^{*}$ are the unique basis
$1$-forms for $V^{*}$ which satisfy
\begin{equation}
\varepsilon_{j}\cdot\varepsilon^{k}=\delta_{j}^{k}. \label{MS.4b}%
\end{equation}

The bases $\{e^{k}\}$ and $\{\varepsilon_{k}\}$ are respectively called the
reciprocal bases of $\{e_{k}\}$ and $\{\varepsilon^{k}\}$ (relatives to the
metric tensor $G$).

Note that $\{\varepsilon_{k}\}$ is the dual basis of $\{e^{k}\},$ i.e.,
\begin{equation}
\varepsilon_{k}(e^{l})=\delta_{k}^{l}, \label{MS.5}%
\end{equation}
an immediate consequence of Eqs.(\ref{MS.4a}) and (\ref{MS.3a}).

From Eqs.(\ref{MS.4a}) and (\ref{MS.4b}), Eqs.(\ref{MS.3a}) and (\ref{MS.3b})
taking into account Eq.(\ref{MS.2}), we easily get that
\begin{align}
\varepsilon_{j}\cdot\varepsilon_{k}  &  =e_{j}\cdot e_{k},\label{MS.5a}\\
e^{j}\cdot e^{k}  &  =G^{jk}=\varepsilon^{j}\cdot\varepsilon^{k}.
\label{MS.5b}%
\end{align}

Using Eq.(\ref{MS.3b}) we get two expansion formulas for $v\in V$
\begin{equation}
v=v\cdot e^{k}e_{k}=v\cdot e_{k}e^{k}. \label{MS.6a}%
\end{equation}

Using Eq.(\ref{MS.4b}) we have that for all $\omega\in V^{\ast}$
\begin{equation}
\omega=\omega\cdot\varepsilon_{k}\varepsilon^{k}=\omega\cdot\varepsilon
^{k}\varepsilon_{k}. \label{MS.6b}%
\end{equation}

Let us take $X\in\bigwedge^{k}V$ with $k\geq2.$ By following analogous steps
to those which allowed us to get Eq.(\ref{EP.5b}) we can now obtain another
expansion formula for $k$-vectors, namely
\begin{equation}
X=\frac{1}{k!}X_{j_{1}\ldots j_{k}}e^{j_{1}}\wedge\ldots\wedge e^{j_{k}},
\label{MS.7}%
\end{equation}
where $X_{j_{1}\ldots j_{k}}=X(\varepsilon_{j_{1}},\ldots,\varepsilon_{j_{k}%
})$ are the so-called $j_{1}\ldots j_{k}$-\emph{covariant components of} $X$
(with respect to the $k$-tensor basis $\{e^{j_{1}}\otimes\ldots\otimes
e^{j_{k}}\}$ with $j_{1},\ldots,j_{k}=1,\ldots,n$).

Next, we will obtain a relation between the $i_{1}\ldots i_{k}$-covariant
components of $X$ and the $j_{1}\ldots j_{k}$-contravariant components of $X.
$ By using Eq.(\ref{MS.6b}) and Eq.(\ref{MS.5a}), a straightforward
calculation yields
\begin{align*}
X(\varepsilon_{i_{1}}\ldots\varepsilon_{i_{k}})  &  =X(e_{i_{1}}\cdot
e_{s_{1}}\varepsilon^{s_{1}},\ldots,e_{i_{k}}\cdot e_{s_{k}}\varepsilon
^{s_{k}})\\
&  =X(\varepsilon^{s_{1}},\ldots,\varepsilon^{s_{k}})(e_{i_{1}}\cdot e_{s_{1}%
})\ldots(e_{i_{k}}\cdot e_{s_{k}})\\
&  =\frac{1}{k!}X(\varepsilon^{j_{1}},\ldots,\varepsilon^{j_{k}})\delta
_{j_{1}\ldots j_{k}}^{s_{1}\ldots s_{k}}(e_{i_{1}}\cdot e_{s_{1}}%
)\ldots(e_{i_{k}}\cdot e_{s_{k}}),
\end{align*}
hence,
\begin{equation}
X_{i_{1}\ldots i_{k}}=\frac{1}{k!}X^{j_{1}\ldots j_{k}}\det\left[
\begin{array}
[c]{ccc}%
e_{i_{1}}\cdot e_{j_{1}} & \ldots & e_{i_{1}}\cdot e_{j_{k}}\\
\ldots & \ldots & \ldots\\
e_{i_{k}}\cdot e_{j_{1}} & \ldots & e_{i_{k}}\cdot e_{j_{k}}%
\end{array}
\right]  . \label{MS.8}%
\end{equation}

\subsection{Scalar Product}

Once a metric structure $(V,G)$ has been given we can equip $\bigwedge^{p}V$
with a scalar product of $p$-vectors. $\bigwedge V$ can then be endowed with a
scalar product of multivectors. This is done as follows.

The scalar product of $X_{p},Y_{p}\in\bigwedge^{p}V,$ namely $X_{p}\cdot
Y_{p}\in\mathbb{R},$ is defined by the axioms:

\textbf{Ax-i} For all $\alpha,\beta\in\mathbb{R}:$%
\begin{equation}
\alpha\cdot\beta=\alpha\beta\text{ (real product of }\alpha\text{ and }%
\beta\text{).} \label{SP.1a}%
\end{equation}

\textbf{Ax-ii} For all $X_{p},Y_{p}\in\bigwedge^{p}V,$ with $p\geq1:$%
\begin{align}
X_{p}\cdot Y_{p}  &  =\frac{1}{p!}X_{p}(\varepsilon^{i_{1}},\ldots
,\varepsilon^{i_{p}})Y_{p}(\varepsilon_{i_{1}},\ldots,\varepsilon_{i_{p}%
}),\nonumber\\
&  =\frac{1}{p!}X_{p}(\varepsilon_{i_{1}},\ldots,\varepsilon_{i_{p}}%
)Y_{p}(\varepsilon^{i_{1}},\ldots,\varepsilon^{i_{p}}), \label{SP.1b}%
\end{align}
where $\{\varepsilon_{i}\}$ is the reciprocal basis of $\{\varepsilon^{i}\},$
as defined by Eq.(\ref{MS.4a}).

It is not difficult to realize that the real number defined by Eq.(\ref{SP.1b}%
) does not depend on the bases $\{\varepsilon_{i}\}$ and $\{\varepsilon^{i}\}$
for calculating it. Indeed, by using the expansion formulas for $1$-forms:
$\omega=\omega\cdot\varepsilon_{j}\varepsilon^{j}$ (relative to any pair of
reciprocal bases $\{\varepsilon^{j}\}$ and $\{\varepsilon_{j}\}$), and
$\omega=\omega\cdot\varepsilon^{i\prime}\varepsilon_{i}^{\prime}$ (relative to
any pair of reciprocal bases $\{\varepsilon^{i\prime}\}$ and $\{\varepsilon
_{i}^{\prime}\}$), we have that
\begin{align*}
&  X_{p}(\varepsilon^{i_{1}\prime},\ldots,\varepsilon^{i_{p}\prime}%
)Y_{p}(\varepsilon_{i_{1}}^{\prime},\ldots,\varepsilon_{i_{p}}^{\prime})\\
&  =X_{p}(\varepsilon^{i_{1}\prime}\cdot\varepsilon_{j_{1}}\varepsilon^{j_{1}%
},\ldots,\varepsilon^{i_{p}\prime}\cdot\varepsilon_{j_{p}}\varepsilon^{j_{p}%
})Y_{p}(\varepsilon_{i_{1}}^{\prime},\ldots,\varepsilon_{i_{p}}^{\prime})\\
&  =X_{p}(\varepsilon^{j_{1}},\ldots,\varepsilon^{j_{p}})Y_{p}(\varepsilon
_{j_{1}}\cdot\varepsilon^{i_{1}\prime}\varepsilon_{i_{1}}^{\prime}%
,\ldots,\varepsilon_{j_{p}}\cdot\varepsilon^{i_{p}\prime}\varepsilon_{i_{p}%
}^{\prime})\\
&  =X_{p}(\varepsilon^{j_{1}},\ldots,\varepsilon^{j_{p}})Y_{p}(\varepsilon
_{j_{1}},\ldots,\varepsilon_{j_{p}}).
\end{align*}

It is a well-defined scalar product on $\bigwedge^{p}V,$ since it is
symmetric, satisfies the distributive laws, has the mixed associativity
property and is non-degenerate i.e., if $X_{p}\cdot Y_{p}=0$ for all $Y_{p},$
then $X_{p}=0.$

We prove here only the non-degeneracy property. Choose, e.g., $Y_{p}\equiv
e^{j_{1}}\wedge\ldots\wedge e^{j_{p}}.$ Thus, all that must proved is that if
$X_{p}\cdot(e^{j_{1}}\wedge\ldots\wedge e^{j_{p}})=0,$ then $X_{p}=0.$

Using Eq.(\ref{EP.5a}), more precisely, its version using the reciprocal basis
vectors $e^{1},\ldots,e^{n},$ the duality condition of $\{e^{k}\}$ and
$\{\varepsilon_{k}\},$ and taking into account Eq.(\ref{EP.2h}), we have that
\begin{align*}
X_{p}\cdot(e^{j_{1}}\wedge\ldots\wedge e^{j_{p}})=0  &  \Rightarrow\frac
{1}{p!}X_{p}(\varepsilon^{i_{1}},\ldots,\varepsilon^{i_{p}})\delta
_{s_{1}\ldots s_{p}}^{j_{1}\ldots j_{p}}\varepsilon_{i_{1}}(e^{s_{1}}%
)\ldots\varepsilon_{i_{p}}(e^{s_{p}})=0\\
&  \Rightarrow\frac{1}{p!}X_{p}(\varepsilon^{i_{1}},\ldots,\varepsilon^{i_{p}%
})\delta_{s_{1}\ldots s_{p}}^{j_{1}\ldots j_{p}}\delta_{i_{1}}^{s_{1}}%
\ldots\delta_{i_{p}}^{s_{p}}=0\\
&  \Rightarrow X_{p}(\varepsilon^{j_{1}},\ldots,\varepsilon^{j_{p}})=0,
\end{align*}
i.e., $X_{p}=0.$

For the special case of vectors Eq.(\ref{SP.1b}) reduces to
\begin{equation}
v\cdot w=\varepsilon^{i}(v)\varepsilon_{i}(w)=\varepsilon_{i}(v)\varepsilon
^{i}(w), \label{SP.1c}%
\end{equation}
i.e., $G=\varepsilon^{i}\otimes\varepsilon_{i}=\varepsilon_{i}\otimes
\varepsilon^{i}.$

Note that Eq.(\ref{SP.1c}) is consistent with Eq.(\ref{MS.1c}). We have indeed
that
\begin{align*}
\varepsilon^{i}(v)\varepsilon_{i}(w)  &  =\varepsilon^{i}(v)G_{ij}%
\varepsilon^{j}(w)=\varepsilon^{i}(v)G(e_{i},e_{j})\varepsilon^{j}(w)\\
&  =G(\varepsilon^{i}(v)e_{i},\varepsilon^{j}(w)e_{j})=G(v,w).
\end{align*}

The well-known formula for the scalar product of simple $k$-vectors can be
easily deduced from Eq.(\ref{SP.1b}). It is:
\begin{equation}
(v_{1}\wedge\ldots\wedge v_{k})\cdot(w_{1}\wedge\ldots\wedge w_{k}%
)=\det\left[
\begin{array}
[c]{ccc}%
v_{1}\cdot w_{1} & \ldots & v_{1}\cdot w_{k}\\
\ldots & \ldots & \ldots\\
v_{k}\cdot w_{1} & \ldots & v_{k}\cdot w_{k}%
\end{array}
\right]  . \label{SP.2}%
\end{equation}

\textbf{Proof}

By using Eq.(\ref{EP.4}) and Eq.(\ref{SP.1c}), and recalling the $k\times k$
determinant formula, $\det\left[  a_{ij}\right]  =\dfrac{1}{k!}\epsilon
^{i_{1}\ldots i_{k}}\epsilon^{j_{1}\ldots j_{k}}a_{i_{1}j_{1}}\ldots
a_{i_{k}j_{k}}.$ A straightforward calculation gives
\begin{align*}
&  (v_{1}\wedge\ldots\wedge v_{k})\cdot(w_{1}\wedge\ldots\wedge w_{k})\\
&  =\dfrac{1}{k!}v_{1}\wedge\ldots\wedge v_{k}(\varepsilon^{s_{1}}%
,\ldots,\varepsilon^{s_{k}})w_{1}\wedge\ldots\wedge w_{k}(\varepsilon_{s_{1}%
},\ldots,\varepsilon_{s_{k}})\\
&  =\dfrac{1}{k!}\epsilon^{i_{1}\ldots i_{k}}\varepsilon^{s_{1}}(v_{i_{1}%
})\ldots\varepsilon^{s_{k}}(v_{i_{k}})\epsilon^{j_{1}\ldots j_{k}}%
\varepsilon_{s_{1}}(w_{j_{1}})\ldots\varepsilon_{s_{k}}(w_{j_{k}})\\
&  =\dfrac{1}{k!}\epsilon^{i_{1}\ldots i_{k}}\epsilon^{j_{1}\ldots j_{k}%
}\varepsilon^{s_{1}}(v_{i_{1}})\varepsilon_{s_{1}}(w_{j_{1}})\ldots
\varepsilon^{s_{k}}(v_{i_{k}})\varepsilon_{s_{k}}(w_{j_{k}})\\
&  =\dfrac{1}{k!}\epsilon^{i_{1}\ldots i_{k}}\epsilon^{j_{1}\ldots j_{k}%
}(v_{i_{1}}\cdot w_{j_{1}})\ldots(v_{i_{k}}\cdot w_{j_{k}}),\\
&  =\det\left[  v_{i}\cdot w_{j}\right]  .\mathbf{\blacksquare}%
\end{align*}

Now, we can generalize Eq.(\ref{MS.6a}) in order to get the expected expansion
formulas for $k$-vectors. For all $X\in\bigwedge^{k}V$ it holds two expansion
formulas
\begin{equation}
X=\frac{1}{k!}X\cdot(e^{j_{1}}\wedge\ldots e^{j_{k}})(e_{j_{1}}\wedge\ldots
e_{j_{k}})=\frac{1}{k!}X\cdot(e_{j_{1}}\wedge\ldots e_{j_{k}})(e^{j_{1}}%
\wedge\ldots e^{j_{k}}). \label{SP.2a}%
\end{equation}

\textbf{Proof}

For $X\in\bigwedge^{k}V$ with $k\geq2,$ by recalling Eq.(\ref{EP.5b}) and
Eq.(\ref{MS.7}) there are unique real numbers $X^{i_{1}\ldots i_{k}}$ and
$X_{i_{1}\ldots i_{k}}$ with $i_{1},\ldots,i_{k}=1,\ldots,n$ such that
\[
X=\frac{1}{k!}X^{i_{1}\ldots i_{k}}e_{i_{1}}\wedge\ldots e_{i_{k}}=\frac
{1}{k!}X_{i_{1}\ldots i_{k}}e^{i_{1}}\wedge\ldots e^{i_{k}}.
\]
Indeed, take the scalar products $X\cdot(e^{j_{1}}\wedge\ldots e^{j_{k}})$.
Using Eq.(\ref{SP.2}), Eq.(\ref{MS.3b}) and Eq.(\ref{EP.2h}) we have that
\begin{align*}
X\cdot(e^{j_{1}}\wedge\ldots e^{j_{k}})  &  =\frac{1}{k!}X^{i_{1}\ldots i_{k}%
}(e_{i_{1}}\wedge\ldots e_{i_{k}})\cdot(e^{j_{1}}\wedge\ldots e^{j_{k}})\\
&  =\frac{1}{k!}X^{i_{1}\ldots i_{k}}\det\left[
\begin{array}
[c]{ccc}%
e_{i_{1}}\cdot e^{j_{1}} & \ldots & e_{i_{1}}\cdot e^{j_{k}}\\
\ldots & \ldots & \ldots\\
e_{i_{k}}\cdot e^{j_{1}} & \ldots & e_{i_{k}}\cdot e^{j_{k}}%
\end{array}
\right] \\
&  =\frac{1}{k!}\delta_{i_{1}\ldots i_{k}}^{j_{1}\ldots j_{k}}X^{i_{1}\ldots
i_{k}}=X^{j_{1}\ldots j_{k}},
\end{align*}
i.e., $X^{j_{1}\ldots j_{k}}=X\cdot(e^{j_{1}}\wedge\ldots e^{j_{k}}).$
Analogously, we can prove that $X_{j_{1}\ldots j_{k}}=X\cdot(e_{j_{1}}%
\wedge\ldots e_{j_{k}})$.$\blacksquare$

\subsection{Scalar Product of Multivectors}

The scalar product of $X,Y\in\bigwedge V,$ namely $X\cdot Y\in\mathbb{R},$ is
defined by
\begin{equation}
X\cdot Y=\overset{n}{\underset{k=0}{\sum}}\pi_{k}(X)\cdot\pi_{k}(Y).
\label{SP.3}%
\end{equation}
Note that on the right side there appears the scalar products of $k$-vectors
with $0\leq k\leq n,$ as defined by Eqs..(\ref{SP.1a}) and (\ref{SP.1b}). It
means that if $X=(X_{0},X_{1},\ldots,X_{n})$ and $Y=(Y_{0},Y_{1},\ldots
,Y_{n}),$ then
\begin{equation}
X\cdot Y=\overset{n}{\underset{k=0}{\sum}}X_{k}\cdot Y_{k}. \label{SP.3a}%
\end{equation}

By using Eqs..(\ref{SP.1a}) and (\ref{SP.1b}) we can easily note that
Eq.(\ref{SP.3}) can still be written as
\begin{align}
X\cdot Y  &  =X_{0}Y_{0}+\overset{n}{\underset{k=1}{\sum}}\frac{1}{k!}%
X_{k}(\varepsilon^{i_{1}},\ldots,\varepsilon^{i_{k}})Y_{k}(\varepsilon_{i_{1}%
},\ldots,\varepsilon_{i_{k}})\nonumber\\
&  =X_{0}Y_{0}+\overset{n}{\underset{k=1}{\sum}}\frac{1}{k!}X_{k}%
(\varepsilon_{i_{1}},\ldots,\varepsilon_{i_{k}})Y_{k}(\varepsilon^{i_{1}%
},\ldots,\varepsilon^{i_{k}}). \label{SP.4}%
\end{align}

It is important to observe that the operation defined by Eq.(\ref{SP.3}) is
indeed a well-defined scalar product on $\bigwedge V$, since it is symmetric,
satisfies the distributive laws, has the mixed associative property and is not
degenerate, i.e., if $X\cdot Y=0$ for all $Y,$ then $X=0.$

Let us take $X_{p}\in\bigwedge^{p}V$ and $Y_{q}\in\bigwedge^{q}V.$ By using
Eq.(\ref{SP.3}) and Eq.(\ref{M.15a}) we have that
\begin{align*}
\tau_{p}(X_{p})\cdot\tau_{q}(Y_{q})  &  =\overset{n}{\underset{k=0}{\sum}}%
\pi_{k}\circ\tau_{p}(X_{p})\cdot\pi_{k}\circ\tau_{q}(Y_{q})\\
&  =\overset{n}{\underset{k=0}{\sum}}\left\{
\begin{array}
[c]{cc}%
0_{k}, & k\neq p\\
X_{p}, & k=p
\end{array}
\right.  \cdot\left\{
\begin{array}
[c]{cc}%
0_{k}, & k\neq q\\
Y_{q}, & k=q
\end{array}
\right.  .
\end{align*}

The sum on the right side of this equation is $0=\overset{n}{\underset
{k=0}{\sum}}0_{k}\cdot0_{k}$ unless there exists some $k_{0}$ with $0\leq
k_{0}\leq n$ such that $k_{0}=p$ and $k_{0}=q.$ Whence, we see that for $p\neq
q$ there is no number $k_{0}$ to satisfy the required conditions. But, for
$p=q$ we have that $k_{0}=p=q$ trivially satisfies them. Thus, from the above
equation it follows that
\begin{equation}
\tau_{p}(X)\cdot\tau_{q}(Y_{q})=\left\{
\begin{array}
[c]{cc}%
0, & \text{if }p\neq q\\
X_{p}\cdot Y_{p}, & \text{if }p=q
\end{array}
\right.  . \label{SP.5}%
\end{equation}

\section{Contracted Products}

The left contracted product of $X_{p}\in\bigwedge^{p}V$ and $Y_{q}\in
\bigwedge^{q}V$ with $0\leq p\leq q\leq n,$ namely $X_{p}\lrcorner Y_{q}%
\in\bigwedge^{q-p}V,$ is defined by the following axioms

\textbf{Ax-i} For all $X_{p},Y_{p}\in\bigwedge^{p}V:$%
\begin{equation}
X_{p}\lrcorner Y_{p}=\widetilde{X_{p}}\cdot Y_{p}=X_{p}\cdot\widetilde{Y_{p}}.
\label{CP.1a}%
\end{equation}

\textbf{Ax-ii} For all $X_{p}\in\bigwedge^{p}V$ and $Y_{q}\in\bigwedge^{q}V $
with $p<q:$%
\begin{align}
X_{p}\lrcorner Y_{q}  &  =\frac{1}{(q-p)!}(\widetilde{X_{p}}\wedge e^{i_{1}%
}\wedge\ldots\wedge e^{i_{q-p}})\cdot Y_{q}e_{i_{1}}\wedge\ldots\wedge
e_{i_{q-p}}\nonumber\\
&  =\frac{1}{(q-p)!}(\widetilde{X_{p}}\wedge e_{i_{1}}\wedge\ldots\wedge
e_{i_{q-p}})\cdot Y_{q}e^{i_{1}}\wedge\ldots\wedge e^{i_{q-p}}. \label{CP.1b}%
\end{align}
The right contracted product of $X_{p}\in\bigwedge^{p}V$ and $Y_{q}%
\in\bigwedge^{q}V$ with $n\geq p\geq q\geq0,$ namely $X_{p}\llcorner Y_{q}%
\in\bigwedge^{p-q}V,$ is defined by the following axioms

\textbf{Ax-i} For all $X_{p},Y_{p}\in\bigwedge^{p}V:$%
\begin{equation}
X_{p}\llcorner Y_{p}=\widetilde{X_{p}}\cdot Y_{p}=X_{p}\cdot\widetilde{Y_{p}}.
\label{CP.2a}%
\end{equation}

\textbf{Ax-ii} For all $X_{p}\in\bigwedge^{p}V$ and $Y_{q}\in\bigwedge^{q}V $
with $p>q:$%
\begin{align}
X_{p}\llcorner Y_{q}  &  =\frac{1}{(p-q)!}X_{p}\cdot(e^{i_{1}}\wedge
\ldots\wedge e^{i_{p-q}}\wedge\widetilde{Y_{q}})e_{i_{1}}\wedge\ldots\wedge
e_{i_{p-q}}\nonumber\\
&  =\frac{1}{(p-q)!}X_{p}\cdot(e_{i_{1}}\wedge\ldots\wedge e_{i_{p-q}}%
\wedge\widetilde{Y_{q}})e^{i_{1}}\wedge\ldots\wedge e^{i_{p-q}}. \label{CP.2b}%
\end{align}
It should be noticed that the $(q-p)$-vector defined by Eq.(\ref{CP.1b}) and
the $(p-q)$-vector defined by Eq.(\ref{CP.2b}) do not depend on the choice of
the reciprocal bases $\{e_{i}\}$ and $\{e^{i}\}$ used for calculating them.

Let us take $X_{p}\in\bigwedge^{p}V$ and $Y_{q}\in\bigwedge^{q}V$ with $p\leq
q.$ For all $Z_{q-p}\in\bigwedge^{q-p}V$ the following identity holds
\begin{equation}
(X_{p}\lrcorner Y_{q})\cdot Z_{q-p}=Y_{q}\cdot(\widetilde{X_{p}}\wedge
Z_{q-p}). \label{CP.3}%
\end{equation}
For $p<q$ Eq.(\ref{CP.3}) follows directly from Eq.(\ref{CP.1b}) and
Eq.(\ref{SP.2a}). But, for $p=q$ it trivially follows by taking into account
Eq.(\ref{CP.1a}), etc.

Let us take $X_{p}\in\bigwedge^{p}V$ and $Y_{q}\in\bigwedge^{q}V$ with $p\geq
q.$ For all $Z_{p-q}\in\bigwedge^{p-q}V$ the following identity holds
\begin{equation}
(X_{p}\llcorner Y_{q})\cdot Z_{p-q}=X_{p}\cdot(Z_{p-q}\wedge\widetilde{Y_{q}%
}). \label{CP.4}%
\end{equation}
For $p>q$ Eq.(\ref{CP.4}) follows directly from Eq.(\ref{CP.2b}) and
Eq.(\ref{SP.2a}). For $p=q$ it follows from Eq.(\ref{CP.2a}).

We present now the basic properties of the elementary contracted products of
$p$-vector with $q$-vector.

For any $V_{p},W_{p}\in\bigwedge^{p}V$ and $X_{q},Y_{q}\in\bigwedge^{q}V$ with
$p\leq q$%
\begin{align}
(V_{p}+W_{p})\lrcorner X_{q}  &  =V_{p}\lrcorner X_{q}+W_{p}\lrcorner
X_{q},\nonumber\\
V_{p}\lrcorner(X_{q}+Y_{q})  &  =V_{p}\lrcorner X_{q}+V_{p}\lrcorner
Y_{q}\text{ (distributive laws).} \label{CP.3a}%
\end{align}

For any $V_{p},W_{p}\in\bigwedge^{p}V$ and $X_{q},Y_{q}\in\bigwedge^{q}V$ with
$p\geq q$%
\begin{align}
(V_{p}+W_{p})\llcorner X_{q}  &  =V_{p}\llcorner X_{q}+W_{p}\llcorner
X_{q},\nonumber\\
V_{p}\llcorner(X_{q}+Y_{q})  &  =V_{p}\llcorner X_{q}+V_{p}\llcorner
Y_{q}\text{ (distributive laws).} \label{CP.4a}%
\end{align}
These distributive laws are immediate consequence of the distributive laws for
the exterior product and the scalar product.

For any $X_{p}\in\bigwedge^{p}V$ and $Y_{q}\in\bigwedge^{q}V$ with $p\leq q $%
\begin{equation}
X_{p}\lrcorner Y_{q}=(-1)^{p(q-p)}Y_{q}\llcorner X_{p}. \label{CP.12}%
\end{equation}
Indeed, by using Eq.(\ref{CP.3}), Eq.(\ref{EP.5f}) and Eq.(\ref{CP.4}) we have
that $(X_{p}\lrcorner Y_{q})\cdot Z_{q-p}=Y_{q}\cdot(\widetilde{X_{p}}\wedge
Z_{q-p})=(-1)^{p(q-p)}Y_{q}\cdot(Z_{q-p}\wedge\widetilde{X_{p}})=(-1)^{p(q-p)}%
(Y_{q}\llcorner X_{p})\cdot Z_{q-p},$ hence, by non-degeneracy of scalar
product, the required result follows.

\subsection{Contracted Product of Multivectors}

The left and right contracted products of $X,Y\in\bigwedge V,$ namely
$X\lrcorner Y\in\bigwedge V$ and $X\llcorner Y\in\bigwedge V,$ are defined by
\begin{align}
X\lrcorner Y  &  =\overset{n}{\underset{k=0}{\sum}}\overset{n-k}%
{\underset{j=0}{\sum}}\tau_{k}(\pi_{j}(X)\lrcorner\pi_{k+j}(Y)).
\label{CP.5a}\\
X\llcorner Y  &  =\overset{n}{\underset{k=0}{\sum}}\overset{n-k}%
{\underset{j=0}{\sum}}\tau_{k}(\pi_{k+j}(X)\llcorner\pi_{j}(Y)). \label{CP.5b}%
\end{align}
Note that on the right side of Eq.(\ref{CP.5a}) there appear the left
contracted products of $j$-vectors with $(k+j)$-vectors, as defined by
Eqs..(\ref{CP.1a}) and (\ref{CP.1b}). On the right side of Eq.(\ref{CP.5b}%
)\ there appear the right contracted product of $(k+j)$-vectors with
$j$-vectors, as defined by Eqs..(\ref{CP.2a}) and (\ref{CP.2b}). It means that%

\begin{align}
\pi_{k}(X\lrcorner Y)  &  =\overset{n-k}{\underset{j=0}{\sum}}\pi
_{j}(X)\lrcorner\pi_{k+j}(Y),\label{CP.5c}\\
\pi_{k}(X\llcorner Y)  &  =\overset{n-k}{\underset{j=0}{\sum}}\pi
_{k+j}(X)\llcorner\pi_{j}(Y)\text{ for each }k=0,1,\ldots,n. \label{CP.5d}%
\end{align}
And, if $X=(X_{0},\ldots,X_{k},\ldots,X_{n})$ and $Y=(Y_{0},\ldots
,Y_{k},\ldots,Y_{n}),$ then
\begin{align}
X\lrcorner Y  &  =(\overset{n}{\underset{j=0}{\sum}}X_{j}\lrcorner
Y_{j},\ldots,\overset{n-k}{\underset{j=0}{\sum}}X_{j}\lrcorner Y_{k+j}%
,\ldots,X_{0}Y_{n}),\label{CP.5e}\\
X\llcorner Y  &  =(\overset{n}{\underset{j=0}{\sum}}X_{j}\llcorner
Y_{j},\ldots,\overset{n-k}{\underset{j=0}{\sum}}X_{k+j}\llcorner Y_{j}%
,\ldots,X_{n}Y_{0}). \label{CP.5f}%
\end{align}

These contracted products are internal laws on $\bigwedge V.$ Each of
$\lrcorner$ and $\llcorner$ satisfies the distributive laws, as easily verify
by recalling Eq.(\ref{CP.3a}) and Eq.(\ref{CP.4a}), but both interior products
are non associative. $\bigwedge V$ endowed with each one of these contracted
products is a non-associative algebra which will be called a \emph{interior
algebra of multivectors}.

Let us take $X_{p}\in\bigwedge^{p}V$ and $Y_{q}\in\bigwedge^{q}V.$ By using
Eq.(\ref{CP.5c}) and Eq.(\ref{M.15a}) we have that
\begin{align*}
\pi_{k}(\tau_{p}(X_{p})\lrcorner\tau_{q}(Y_{q}))  &  =\overset{n-k}%
{\underset{j=0}{\sum}}\pi_{j}\circ\tau_{p}(X_{p})\lrcorner\pi_{k+j}\circ
\tau_{q}(Y_{q})\\
&  =\overset{n-k}{\underset{j=0}{\sum}}\left\{
\begin{array}
[c]{cc}%
0_{j}, & j\neq p\\
X_{p}, & j=p
\end{array}
\right.  \lrcorner\left\{
\begin{array}
[c]{cc}%
0_{k+j}, & k+j\neq q\\
Y_{q}, & k+j=q
\end{array}
\right.  .
\end{align*}

The sum on the right side of the equation above is $0_{k}=\overset
{n-k}{\underset{j=0}{%
{\displaystyle\sum}
}}0_{j}\lrcorner0_{k+j},$ for each $k=0,1,\ldots,n$ unless there exists some
$k_{0}$ with $0\leq k_{0}\leq n$ such that $k_{0}+p=q.$ Whence, we see that
for $p>q$ there is no number $k_{0}$ which can satisfy the required
conditions. But, for $p\leq q$ we have that $k_{0}=q-p$ is the unique number
to satisfy them. Then, from the equation above it follows that
\begin{align*}
\pi_{k}(\tau_{p}(X_{p})\lrcorner\tau_{q}(Y_{q}))  &  =0_{k},\text{ for }p>q\\
\pi_{q-p}(\tau_{p}(X_{p})\lrcorner\tau_{q}(Y_{q}))  &  =X_{p}\lrcorner
Y_{q},\text{ for }p\leq q.
\end{align*}
Hence, by taking into account Eq.(\ref{M.8a}) we finally obtain
\begin{equation}
\tau_{p}(X_{p})\lrcorner\tau_{q}(Y_{q})=\left\{
\begin{array}
[c]{cc}%
0, & \text{if }p>q\\
\tau_{q-p}(X_{p}\lrcorner Y_{q}), & \text{if }p\leq q
\end{array}
\right.  . \label{CP.6}%
\end{equation}

By following analogous steps to those which allowed us to arrive to
Eq.(\ref{CP.6}) we can also obtain
\begin{equation}
\tau_{p}(X_{p})\llcorner\tau_{q}(Y_{q})=\left\{
\begin{array}
[c]{cc}%
0, & \text{if }p<q\\
\tau_{p-q}(X_{p}\llcorner Y_{q}), & \text{if }p\geq q
\end{array}
\right.  . \label{CP.6a}%
\end{equation}

We finalize this section by presenting two noticeable formulas involving the
contracted products and the scalar product, and two other remarkable formulas
relating the contracted products to the exterior product and scalar product.

For any $X,Y,Z\in\bigwedge V$
\begin{align}
(X\lrcorner Y)\cdot Z  &  =Y\cdot(\widetilde{X}\wedge Z),\label{CP.8a}\\
(X\llcorner Y)\cdot Z  &  =X\cdot(Z\wedge\widetilde{Y}). \label{CP.8b}%
\end{align}

\textbf{Proof}

We only give the proof for the first statement, the other being analogous.
Form the summation identity $\overset{n}{\underset{k=0}{\sum}}\overset
{n-k}{\underset{j=0}{\sum}}A_{k+j}\cdot(B_{j}\wedge C_{k})=\overset
{n}{\underset{k=0}{\sum}}\overset{k}{\underset{j=0}{\sum}}A_{k}\cdot
(B_{j}\wedge C_{k-j})$ and by using Eqs.(\ref{SP.3}), (\ref{CP.5c}),
(\ref{CP.3}) and (\ref{EP.6a}), a straightforward calculation yields
\begin{align*}
(X\lrcorner Y)\cdot Z  &  =\overset{n}{\underset{k=0}{\sum}}\pi_{k}(X\lrcorner
Y)\cdot Z_{k}=\overset{n}{\underset{k=0}{\sum}}(\overset{n-k}{\underset
{j=0}{\sum}}X_{j}\lrcorner Y_{k+j})\cdot Z_{k}\\
&  =\overset{n}{\underset{k=0}{\sum}}\overset{n-k}{\underset{j=0}{\sum}}%
(X_{j}\lrcorner Y_{k+j})\cdot Z_{k}=\overset{n}{\underset{k=0}{\sum}}%
\overset{n-k}{\underset{j=0}{\sum}}Y_{k+j}\cdot(\widetilde{X_{j}}\wedge
Z_{k})\\
&  =\overset{n}{\underset{k=0}{\sum}}\overset{k}{\underset{j=0}{\sum}}%
Y_{k}\cdot(\widetilde{X_{j}}\wedge Z_{k-j})=\overset{n}{\underset{k=0}{\sum}%
}Y_{k}\cdot(\overset{k}{\underset{j=0}{\sum}}\widetilde{X_{j}}\wedge
Z_{k-j})\\
&  =\overset{n}{\underset{k=0}{\sum}}Y_{k}\cdot\pi_{k}(\widetilde{X}\wedge
Z)=Y\cdot(\widetilde{X}\wedge Z).\blacksquare
\end{align*}

For any $X,Y,Z\in\bigwedge V$
\begin{align}
X\lrcorner(Y\lrcorner Z)  &  =(X\wedge Y)\cdot Z,\label{CP.9a}\\
(X\llcorner Y)\llcorner Z  &  =X\cdot(Y\wedge Z). \label{CP.9b}%
\end{align}

\textbf{Proof}

We prove only Eq.(\ref{CP.9a}). The proof of Eq.(\ref{CP.9b}) is analogous and
will be left to the reader.

Let $W\in\bigwedge V$. By using Eq.(\ref{CP.8a}) and the associative law for
the exterior product of multivectors, we have that
\begin{align*}
(X\lrcorner(Y\lrcorner Z))\cdot W  &  =(Y\lrcorner Z)\cdot(\widetilde{X}\wedge
W)=Z\cdot((\widetilde{Y}\wedge\widetilde{X})\wedge W)\\
&  =Z\cdot(\widetilde{(X\wedge Y)}\wedge W)=((X\wedge Y)\lrcorner Z)\cdot
W.\blacksquare
\end{align*}
Hence, by the non-degeneracy of the $G$-scalar product, it follows the
required result.

\section{Clifford Product and $\mathcal{C}\ell(V,G)$}

The two interior algebras together with the exterior algebra allow us to
define a \emph{Clifford product} of multivectors which is also an internal law
on $\bigwedge V.$ The Clifford product of $X,Y\in\bigwedge V,$ denoted by
juxtaposition $XY\in\bigwedge V,$ is defined by the following axioms

\textbf{Ax-ci} For all $\alpha\in\mathbb{R}$ and $X\in\bigwedge V:$%
\begin{equation}
\tau_{0}(\alpha)X=\alpha X\text{ (scalar multiplication of }X\text{ by }%
\alpha\text{).} \label{CL.1a}%
\end{equation}

\textbf{Ax-cii} For all $v\in V$ and $X\in\bigwedge V:$%
\begin{align}
\tau_{1}(v)X  &  =\tau_{1}(v)\lrcorner X+\tau_{1}(v)\wedge X,\label{CL.1b}\\
X\tau_{1}(v)  &  =X\llcorner\tau_{1}(v)+X\wedge\tau_{1}(v). \label{CL.1c}%
\end{align}

\textbf{Ax-ciii} For all $X,Y,Z\in\bigwedge V:$%
\begin{equation}
(XY)Z=X(YZ). \label{CL.1d}%
\end{equation}

The Clifford product is distributive and associative. $\bigwedge V$ endowed
with this Clifford product is an associative algebra which will be called
the\emph{\ geometric algebra of multivectors} associated to a metric structure
$(V,G).$ It will be denoted by $\mathcal{C}\ell(V,G).$

Now, due to the existence of a linear isomorphism between the $k$-vectors and
the $k$-homogeneous multivector as given by Eqs..(\ref{M.6a}) and (\ref{M.6b},
and by recalling the remarkable propositions about the exterior product, the
scalar product and the contracted products of multivectors as given by
Eq.(\ref{EP.7}), Eq.(\ref{SP.5}) and Eqs..(\ref{CP.6}) and (\ref{CP.6a}), we
can introduce a notational convention which are more convenient for performing
calculations with multivectors using the geometric algebra.

In whichever addition or products of multivectors all $k$-homogeneous
multivectors $\tau_{k}(X_{k})$ will be \textit{identified} with its
corresponding $k$-vector $X_{k}.$

\textbf{i.} Any addition of multivectors $\tau_{p}(X_{p})+\tau_{q}(Y_{q}),$
$\tau_{p}(X_{p})+Y$ and $X+\tau_{q}(Y_{q})$ will be respectively denoted by
$X_{p}+Y_{q},$ $X_{p}+Y$ and $X+Y_{q}.$

\textbf{ii.} Any product of multivectors $\tau_{p}(X_{p})*\tau_{q}(Y_{q}),$
$\tau_{p}(X_{p})*Y$ and $X*\tau_{q}(Y_{q}),$ where $*$ means either
$(\wedge),$ $(\cdot),$ $(\lrcorner,\llcorner)$ or $($\emph{Clifford
product}$),$ will be respectively denoted by $X_{p}*Y_{q},$ $X_{p}*Y$ and
$X*Y_{q}.$

\section{Euclidean and pseudo-Euclidean Geometric Algebras}

Let us equip $V$ with an arbitrary (but fixed once for all) euclidean metric
$G_{E}$, i.e., a metric tensor on $V$ with the strong condition of being
positive definite, i.e.,
\begin{equation}
G_{E}(v,v)\geq0\text{ for all }v\in V\text{ and if }G_{E}(v,v)=0,\text{ then
}v=0. \label{GA.1}%
\end{equation}

$V$ endowed with an euclidean metric $G_{E},$ i.e., $(V,G_{E}),$ is called an
euclidean metric structure for $V.$ Sometimes, $(V,G_{E})$ is said to be an
euclidean space.

Associated to $(V,G_{E})$ an euclidean scalar product of vectors $v,w\in V$ is
given by
\begin{equation}
v\underset{G_{E}}{\cdot}w=G_{E}(v,w). \label{GA.2}%
\end{equation}
We introduce also an euclidean scalar product of $p$-vectors $X_{p},Y_{p}%
\in\bigwedge^{p}V$ and euclidean scalar product of multivectors $X,Y\in
\bigwedge V,$ namely $X_{p}\underset{G_{E}}{\cdot}Y_{p}\in\mathbb{R}$ and
$X\underset{G_{E}}{\cdot}Y\in\mathbb{R},$ using respectively the
Eqs..(\ref{SP.1a}) and (\ref{SP.1b}), and Eq.(\ref{SP.3}).

Let us take any metric tensor $G$ on the vector space $V.$ Associated to the
metric structure $(V,G)$ a scalar product of vectors $v,w\in V$ is represented
by
\begin{equation}
v\underset{G}{\cdot}w=G(v,w). \label{GA.3}%
\end{equation}
Of course, the corresponding scalar product of $p$-vectors $X_{p},Y_{p}%
\in\bigwedge^{p}V$ and scalar product of multivectors $X,Y\in\bigwedge V,$
namely $X_{p}\underset{G}{\cdot}Y_{p}\in\mathbb{R}$ and $X\underset{G}{\cdot
}Y\in\mathbb{R},$ are defined respectively by Eqs.(\ref{SP.1a}) and
(\ref{SP.1b}) and Eq.(\ref{SP.3}).

We will find a relationship between $(V,G)$ and $(V,G_{E})$, thereby showing
how an arbitrary $G$-scalar product on $\bigwedge^{p}V$ and $\bigwedge V$ is
related to a $G_{E}$-scalar products on $\bigwedge^{p}V$ and $\bigwedge V$.

Choose once and for all a fundamental euclidean metric structure $(V,G_{E})$.
For any metric tensor $G$ there exists an unique linear operator $g$ such that
for all $v,w\in V$
\begin{equation}
v\underset{G}{\cdot}w=g(v)\underset{G_{E}}{\cdot}w. \label{GA.4}%
\end{equation}
Such $g$ is given by
\begin{equation}
g(v)=(v\underset{G}{\cdot}e_{k})\underset{G_{E}}{e^{k}}=(v\underset{G}{\cdot
}\underset{G_{E}}{e^{k}})e_{k}, \label{GA.5}%
\end{equation}
where $\{e_{k}\}$ is any basis of $V,$ and $\{\underset{G_{E}}{e^{k}}\}$ is
its reciprocal basis with respect to $(V,G_{E}),$ i.e., $e_{k}\underset{G_{E}%
}{\cdot}\underset{G_{E}}{e^{l}}=\delta_{k}^{l}.$ Note that the vector $g(v) $
does not depend on the basis $\{e_{k}\}$ chosen for calculating it.

We prove now that $g(v)$ as given by Eq.(\ref{GA.5}) satisfies Eq.(\ref{GA.4}%
). Using Eq.(\ref{MS.6a}) we have,
\[
g(v)\underset{G_{E}}{\cdot}w=(v\underset{G}{\cdot}e_{k})(\underset{G_{E}%
}{e^{k}}\underset{G_{E}}{\cdot}w)=(v\underset{G}{\cdot}(\underset{G_{E}}%
{e^{k}}\underset{G_{E}}{\cdot}w)e_{k})=v\underset{G}{\cdot}w.
\]

Now, suppose that there is some $g^{\prime}$ which satisfies Eq.(\ref{GA.4}),
i.e., $v\underset{G}{\cdot}w=g^{\prime}(v)\underset{G_{E}}{\cdot}w.$ Then, by
using once again Eq.(\ref{MS.6a}) we have that
\[
g^{\prime}(v)=(g^{\prime}(v)\underset{G_{E}}{\cdot}e_{k})\underset{G_{E}%
}{e^{k}}=(v\underset{G}{\cdot}e_{k})\underset{G_{E}}{e^{k}}=g(v),
\]
i.e., $g^{\prime}=g.$

So the existence and the uniqueness of such a linear operator $g$ are proved.

Since $G$ is a symmetric covariant $2$-tensor over $V$, i.e., $G(v,w)=G(w,v)$
$\forall v,w\in V$, it follows from Eq.(\ref{GA.4}) that $g$ is an
\emph{adjoint symmetric} linear operator with respect to $(V,G_{E})$, i.e.,
\begin{equation}
g(v)\underset{G_{E}}{\cdot}w=v\underset{G_{E}}{\cdot}g(w), \label{GA.6}%
\end{equation}
The property expressed by Eq.(\ref{GA.6}) is coded by the equation
$g=g^{\dagger(G_{E})}.$

Since $G$ is a non-degenerate covariant $2$-tensor over $V$ (i.e., if
$G(v,w)=0$ $\forall w\in V,$ then $v=0$) it follows that $g$ is a non-singular
(invertible) linear operator. Its inverse linear operator is given by the
noticeable formula
\begin{equation}
g^{-1}(v)=G^{jk}(v\underset{G_{E}}{\cdot}e_{j})e_{k}, \label{GA.7}%
\end{equation}
where $G^{jk}$ are the $jk$-entries of the inverse matrix of $\left[
G_{jk}\right]  $ with $G_{jk}\equiv G(e_{j},e_{k}).$ Note that the vector
$g^{-1}(v) $ does not depend on the basis $\{e_{k}\}$ chosen for its calculation.

We must prove that indeed $g^{-1}\circ g=g\circ g^{-1}=i_{V},$ where $i_{V}$
is the identity function for $V.$

By using Eq.(\ref{GA.7}), Eq.(\ref{GA.4}), Eq.(\ref{MS.3a}) for $(V,G)$ and
Eq.(\ref{MS.6a}) for $(V,G),$ we have that
\[
g^{-1}\circ g(v)=G^{jk}(g(v)\underset{G_{E}}{\cdot}e_{j})e_{k}=G^{jk}%
(v\underset{G}{\cdot}e_{j})e_{k}=(v\underset{G}{\cdot}e_{j})\underset{G}%
{e^{k}}=v,
\]
i.e., $g^{-1}\circ g=i_{V}.$

By using Eq.(\ref{GA.7}), Eq.(\ref{MS.6a}) for $(V,G_{E}),$ Eq.(\ref{GA.4}),
Eq.(\ref{MS.3a}) for $(V,G),$ Eq.(\ref{MS.3b}) for $(V,G)$ and Eq.(\ref{MS.6a}%
) for $(V,G_{E}),$ we have that
\begin{align*}
g\circ g^{-1}(v)  &  =G^{jk}(v\underset{G_{E}}{\cdot}e_{j})g(e_{k}%
)=G^{jk}(v\underset{G_{E}}{\cdot}e_{j})(g(e_{k})\underset{G_{E}}{\cdot}%
e_{l})\underset{G_{E}}{e^{l}}\\
&  =G^{jk}(v\underset{G_{E}}{\cdot}e_{j})(e_{k}\underset{G}{\cdot}%
e_{l})\underset{G_{E}}{e^{l}}=(v\underset{G_{E}}{\cdot}e_{j})(\underset
{G}{e^{j}}\underset{G}{\cdot}e_{l})\underset{G_{E}}{e^{l}}\\
&  =(v\underset{G_{E}}{\cdot}e_{j})\delta_{l}^{j}\underset{G_{E}}{e^{l}}=v,
\end{align*}
i.e., $g\circ g^{-1}=i_{V}.$

It should be remarked that such $g$ only depends on the choice of the
fundamental euclidean structure $(V,G_{E}).$ However, $g$ codifies all the
geometric information contained in $G.$ Such $g$ will be called the
\emph{metric operator} for $G.$

Now, we show how the scalar product $X_{p}\underset{G}{\cdot}Y_{p}$ is related
to the scalar product $X_{p}\underset{G_{E}}{\cdot}Y_{p}.$

For any simple $k$-vectors $v_{1}\wedge\ldots\wedge v_{k}\in\bigwedge^{k}V $
and $w_{1}\wedge\ldots\wedge w_{k}\in\bigwedge^{k}V$ it holds
\begin{equation}
(v_{1}\wedge\ldots v_{k})\underset{G}{\cdot}(w_{1}\wedge\ldots w_{k}%
)=\underline{g}(v_{1}\wedge\ldots v_{k})\underset{G_{E}}{\cdot}(w_{1}%
\wedge\ldots w_{k}), \label{GA.8}%
\end{equation}
where $\underline{g}$ is the so-called outermorphism (or exterior
power)\footnote{$\underline{g}$ is the unique linear operator on $\bigwedge V$
which satisfies the following conditions: (\textbf{i}) $\alpha\in
\mathbb{R}:\underline{g}(\alpha)=\alpha,$ (\textbf{ii}) $v\in V:\underline
{g}(v)=g(v) $ and (\textbf{iii}) $X,Y\in\bigwedge V:\underline{g}(X\wedge
Y)=\underline{g}(X)\wedge\underline{g}(Y). $} of $g.$

\textbf{Proof}

We will use Eq.(\ref{SP.2}) for $(V,G)$ and $(V,G_{E}).$ By using
Eq.(\ref{GA.4}) and a fundamental property of any outermorphism, a
straightforward calculation yields
\begin{align*}
(v_{1}\wedge\ldots v_{k})\underset{G}{\cdot}(w_{1}\wedge\ldots w_{k})  &
=\det\left[
\begin{array}
[c]{ccc}%
v_{1}\underset{G}{\cdot}w_{1} & \ldots & v_{1}\underset{G}{\cdot}w_{k}\\
\ldots & \ldots & \ldots\\
v_{k}\underset{G}{\cdot}w_{1} & \ldots & v_{k}\underset{G}{\cdot}w_{k}%
\end{array}
\right] \\
&  =\det\left[
\begin{array}
[c]{ccc}%
g(v_{1})\underset{G_{E}}{\cdot}w_{1} & \ldots & g(v_{1})\underset{G_{E}}%
{\cdot}w_{k}\\
\ldots & \ldots & \ldots\\
g(v_{k})\underset{G_{E}}{\cdot}w_{1} & \ldots & g(v_{k})\underset{G_{E}}%
{\cdot}w_{k}%
\end{array}
\right] \\
&  =(g(v_{1})\wedge\ldots g(v_{k}))\underset{G_{E}}{\cdot}(w_{1}\wedge\ldots
w_{k})\\
&  =\underline{g}(v_{1}\wedge\ldots v_{k})\underset{G_{E}}{\cdot}(w_{1}%
\wedge\ldots w_{k}).\blacksquare
\end{align*}

For any $k$-vectors $X_{k},Y_{k}\in\bigwedge^{k}V$ it holds
\begin{equation}
X_{k}\underset{G}{\cdot}Y_{k}=\underline{g}(X_{k})\underset{G_{E}}{\cdot}%
Y_{k}. \label{GA.9}%
\end{equation}

\textbf{Proof}

From the distributive laws and the mixed associative property for the scalar
product of $k$-vectors with $k\geq1$, we get using Eq.(\ref{EP.5b}) and
Eq.(\ref{GA.8}) that
\begin{align*}
X_{k}\underset{G}{\cdot}Y_{k}  &  =(\frac{1}{k!})^{2}X_{k}^{i_{1}\ldots i_{k}%
}Y_{k}^{j_{1}\ldots j_{k}}(e_{i_{1}}\wedge\ldots e_{i_{k}})\underset{G}{\cdot
}(e_{j_{1}}\wedge\ldots e_{j_{k}})\\
&  =(\frac{1}{k!})^{2}X_{k}^{i_{1}\ldots i_{k}}Y_{k}^{j_{1}\ldots j_{k}%
}\underline{g}(e_{i_{1}}\wedge\ldots e_{i_{k}})\underset{G_{E}}{\cdot
}(e_{j_{1}}\wedge\ldots e_{j_{k}})\\
&  =\underline{g}(\frac{1}{k!}X_{k}^{i_{1}\ldots i_{k}}e_{i_{1}}\wedge\ldots
e_{i_{k}})\underset{G_{E}}{\cdot}(\frac{1}{k!}Y_{k}^{j_{1}\ldots j_{k}%
}e_{j_{1}}\wedge\ldots e_{j_{k}})\\
&  =\underline{g}(X_{k})\underset{G_{E}}{\cdot}Y_{k}.
\end{align*}
Note that for the special case of the scalar product of scalars,
Eq.(\ref{GA.9}) trivially holds by using Eq.(\ref{SP.1a}) and taking into
account a fundamental property of any outermorphism.$\blacksquare$

Next, we show how the scalar product $X\underset{G}{\cdot}Y$ is related to the
scalar product $X\underset{G_{E}}{\cdot}Y.$

For any multivectors $X,Y\in\bigwedge V$ it holds
\begin{equation}
X\underset{G}{\cdot}Y=\underline{g}(X)\underset{G_{E}}{\cdot}Y. \label{GA.10}%
\end{equation}

\textbf{Proof}

We will use Eq.(\ref{SP.3}) for $(V,G)$ and $(V,G_{E}).$ By using
Eq.(\ref{GA.9}) and taking into account that $\underline{g}\circ\pi_{k}%
=\pi_{k}\circ\underline{g}$ (the grade-preserving property for the
outermorphism), a straightforward calculation yields
\begin{align*}
X\underset{G}{\cdot}Y  &  =\overset{n}{\underset{k=0}{\sum}}\pi_{k}%
(X)\underset{G}{\cdot}\pi_{k}(Y)=\overset{n}{\underset{k=0}{\sum}}%
\underline{g}\circ\pi_{k}(X)\underset{G_{E}}{\cdot}\pi_{k}(Y)\\
&  =\overset{n}{\underset{k=0}{\sum}}\pi_{k}\circ\underline{g}(X)\underset
{G_{E}}{\cdot}\pi_{k}(Y)=\underline{g}(X)\underset{G_{E}}{\cdot}Y.\blacksquare
\end{align*}

The $G$-contracted products are related to the $G_{E}$-contracted products by
two noticeable formulas.

For any $X,Y\in\bigwedge V$
\begin{align}
X\underset{G}{\lrcorner}Y  &  =\underline{g}(X)\underset{G_{E}}{\lrcorner
}Y,\label{GA.11a}\\
X\underset{G}{\llcorner}Y  &  =X\underset{G_{E}}{\llcorner}\underline{g}(Y).
\label{GA.11b}%
\end{align}

\textbf{Proof}

We only give the proof of the first statement, the other is analogous. By
using Eq.(\ref{CP.8a}) for the metric structure $(V,G)$ we have
that\footnote{Recall that $\underline{g}^{-1}\equiv\underline{(g^{-1}%
)}=(\underline{g})^{-1}.$} for all $X,Y,Z\in\bigwedge V$
\[
(X\underset{G}{\lrcorner}Y)\underset{G}{\cdot}\underline{g}^{-1}%
(Z)=Y\underset{G}{\cdot}(\widetilde{X}\wedge\underline{g}^{-1}(Z)),
\]
by using Eq.(\ref{GA.10})
\[
\underline{g}(X\underset{G}{\lrcorner}Y)\underset{G_{E}}{\cdot}\underline
{g}^{-1}(Z)=\underline{g}(Y)\underset{G_{E}}{\cdot}(\widetilde{X}%
\wedge\underline{g}^{-1}(Z)).
\]
Since the outermorphism of an adjoint symmetric linear operator on $V$ is
itself an adjoint symmetric linear operator on $\bigwedge V$, then recalling a
fundamental property of the outermorphism, it follows that
\begin{align*}
(X\underset{G}{\lrcorner}Y)\underset{G_{E}}{\cdot}Z  &  =Y\underset{G_{E}%
}{\cdot}(\underline{g}(\widetilde{X})\wedge Z),\\
&  =Y\underset{G_{E}}{\cdot}(\widetilde{\underline{g}(X)}\wedge Z).
\end{align*}
Now, by using Eq.(\ref{CP.8a}) for the euclidean metric structure $(V,G_{E})$
we have that for all $X,Y,Z\in\bigwedge V$
\[
(X\underset{G}{\lrcorner}Y)\underset{G_{E}}{\cdot}Z=(\underline{g}%
(X)\underset{G_{E}}{\lrcorner}Y)\underset{G_{E}}{\cdot}Z,
\]
hence, by non-degeneracy of the $G_{E}$-scalar product of multivectors, it
follows that
\[
(X\underset{G}{\lrcorner}Y)=\underline{g}(X)\underset{G_{E}}{\lrcorner
}Y.\blacksquare
\]

\section{Conclusions}

In this paper, the first in a series of eight we start setting the algebraic
basis of the geometrical calculus. We introduced the concept of multivectors
in section 2 and their exterior algebra in section 3. In section 4 we
introduced the key concept of the scalar product of multivectors and in
section 5 the concepts of right and left contractions and interior algebras.
In section 6 we give a definition of a general \emph{real} Clifford (or
geometrical) algebra of multivectors and in section 7 we study in details the
euclidean and pseudo-euclidean geometrical algebras. We emphasize that our
presentation of the geometric (or Clifford) algebras has been devised in order
to give to the reader a powerful computational tool and to permit the
introduction of the theory of extensors in a natural and simple way. To better
achieve our objective we gave in this paper and in the following ones many
details of the calculations (tricks of the trade). To readers interested in
aspects of the general theory of Clifford algebras not covered here we
strongly recommend the excellent textbooks \cite{12,13,14,15}.

\textbf{Acknowledgments: } V. V. Fern\'{a}ndez and A. M. Moya are very
grateful to Mrs. Rosa I. Fern\'{a}ndez who gave to them material and spiritual
support at the starting time of their research work. This paper could not have
been written without her inestimable help.


\begin{thebibliography}{99}                                                                                               %


\bibitem {15}Ablamowicz, R., Baylis, W. E., Branson, T., Lounesto, P.,
Porteous, I., Ryan, J., Selig, J.M., Sobczyk, G., in Ablamowicz, R. and
Sobczyk, G. (eds.), \textit{Lectures on Clifford (Geometrical) Algebras and
Applications, }Birkh\"{a}user, Boston (2004).

\bibitem {6}Fern\'{a}ndez, V. V., Moya, A. M., and Rodrigues, W. A., Jr.,
\emph{Derivative Operators in Metric and Geometric Structures}, submitted for publication.

\bibitem {7}Fern\'{a}ndez, V. V., Moya, A. M., Rodrigues, W. A., Jr. and
Rocha, R., \emph{Riemann and Ricci Fields in Geometric Structures}, submitted
for publication.

\bibitem {10}\noindent Fern\'{a}ndez, V. V., Moya, A. M., and Rodrigues, W.
A., Jr., Mutivector and Extensor Calculus, Special Issue of \emph{Adv. in
Appl. Clifford Algebras }\textbf{11}(S3), 1-103\emph{\ }(2001).

\bibitem {8}Hestenes, D., and Sobcyk, G., \emph{Clifford} \emph{Algebras to
Geometrical Calculus}, D. Reidel Publ. Co., Dordrecht, 1984.

\bibitem {9}Lasenby, A., Doran, C., and Gull, S., Gravity, Gauge Theories and
Geometric Algebras, \emph{Phil. Trans. R. Soc. }\textbf{356}, 487-582 (1998).

\bibitem {12}Lounesto, P., \emph{Clifford Algebras and Spinors}, London Math.
Soc., Lecture Notes Series \textbf{239}, Cambridge University Press,
Cambridge, 1997; second edition 2001.

\bibitem {moro}Mosna, R. A. and Rodrigues, W. A. Jr., The Bundles of Algebraic
and Dirac-Hestenes Spinor Fields, \textit{J. Math. Phys.} \textbf{45},
2945-2966 (2004).

\bibitem {1}Moya, A. M., Fern\'{a}ndez, V. V., and Rodrigues, W. A., Jr.,
\emph{Metric and Gauge Extensors}, submited for publication.

\bibitem {2}Moya, A. M., Fern\'{a}ndez, V. V., and Rodrigues, W. A., Jr.,
\emph{Extensors in Geometric Algebra}, submitted to publication.

\bibitem {3}Moya, A. M., Fern\'{a}ndez, V. V., and Rodrigues, W. A., Jr.
\emph{Multivector and Extensor Fields on Arbitrary Manifolds}, submitted for publication.

\bibitem {4}Moya, A. M., Fern\'{a}ndez, V. V., and Rodrigues, W. A., Jr.,
\emph{Covariant Derivatives of Multivector and Extensor Fields and Intrinsic
Cartan Theory}, submitted for publication.

\bibitem {13}Porteous, I. R., \emph{Clifford Algebras and the Classical
Groups}, Cambridge Studies in Advanced Mathematics \textbf{50}, Cambridge
University Press, Cambridge, 1995; second edition 2001.

\bibitem {14}Porteous, I. R., \emph{Topological Geometry}, Van Nostrand
Reinhold, London, 1969; 2nd edition, Cambridge University Press, Cambridge, 1981.

\bibitem {5}Rodrigues, W. A., Jr., Fern\'{a}ndez, V. V., and Moya, A. M.,
\emph{Metric Compatible Covariant Derivatives}, submitted for publication.

\bibitem {rod041}Rodrigues, W. A. Jr., Algebraic and Dirac-Hestenes Spinors
and Spinor Fields, \textit{J. Math. Phys.} \textbf{45}, 2908-2944 (2004).

\bibitem {rodoliv2006}Rodrigues, W. A. Jr., and Oliveira, E. Capelas, The Many
Faces of Maxwell, Dirac and Einstein Equations, RP 56-05 IMECC-UNICAMP, http://www.ime.unicamp.br/rel\_pesq/2005/rp56-05.html

\bibitem {11}Sobczyk, G., Direct Integration, in Baylis, W. E. (ed.),
\emph{Clifford (Geometric) Algebras with Applications in Physics, Mathematics
and Engineering}, pp. 53-64, Birkh\"{a}user, Boston, 1999.
\end{thebibliography}
\end{document}